\documentstyle{amsart}  %
\newsymbol\varkappa 207B
\newcommand{\mathcal}[1]{\cal{#1}}
\newcommand{\mathbb}[1]{\Bbb{#1}}
\newcommand{\mathfrak}[1]{\frak{#1}}

\newcounter{mycomment}
\newcommand\Hilb{\mathcal{H}}

\newcommand\slim{\operatorname{s-}\lim}
\newcommand\Ran{\operatorname{Ran}}

\renewcommand{\Re}{\operatorname{\text{\shape{n}\rm Re}}}

 \newtheorem{lemma}{Lemma}[section]
 \newtheorem{theorem}[lemma]{Theorem} 
 
\newtheorem{proposition}[lemma]{Proposition}

\theoremstyle{definition}

\newtheorem{definition}[lemma]{Definition}

\theoremstyle{remark}

\newtheorem{remark}[lemma]{Remark}

\numberwithin{equation}{section}

\begin{document}

\title[Scattering problems for perturbations of   media]{Trace-class approach in scattering
problems for perturbations of media}

\author{D. R. Yafaev}
\address{Department of Mathematics, University of Rennes -- I, Campus de Beaulieu, Rennes,
 35042 FRANCE}
\email{yafaev@@univ-rennes1.fr}

\maketitle 

\begin{abstract}
 We consider the
operators $H_0=M_0^{-1}(x) P(D)$ and $H =M^{-1} (x) \\ P(D)$ where $M_0 (x)$ and $M  (x)$ are
positively definite bounded matrix-valued functions and $P(D)$ is an
 elliptic differential operator. 
Our main result is that   the wave operators for the pair $H_0$, $H$ exist
and are complete if the difference $ M(x)-M_0(x)=O(|x|^{-\rho})$, $\rho>d$, as $|x|\rightarrow
\infty$. Our point is that no special assumptions  on $M_0(x)$ are required.
Similar results are obtained in scattering theory for the wave equation.
\end{abstract}

 \section{Introduction}
 
  There are two essentially different methods in scattering theory: the smooth and the
trace-class (see \cite{Y2}, for a more thorough discussion). The first of them originated in the Friedrichs-Faddeev model where
a perturbation of the operator of multiplication $H_0$ by an integral operator $V$ with smooth
kernel is considered. The second goes back to the fundamental Kato-Rosenblum theorem which states
that the wave operators for the pair $H_0$, $H=H_0+V$ exist for a perturbation $V$ from the trace 
class. In applications to differential operators the smooth method works if  
the  operator $H_0$ has constant coefficients and coefficients of $V$  tend  to zero sufficiently
rapidly at infinity. The trace method does not require that coefficients of the operator $H_0$
be constant,  but its assumptions on the fall-off of coefficients of $V$  at infinity are more
stringent. The advantages of the trace-class method were discussed in \cite{BY} for the case where
$H_0=-\Delta + v_0(x)$ is the Schr\"odinger operator with an arbitrary bounded potential
$v_0(x)$, $x\in{\Bbb R}^d$, and $V$ is a first-order differential operator with coefficients bounded
by $|x|^{-\rho}$, $\rho>d$, as $|x|\rightarrow \infty$.

 Our goal here is to apply
the trace-class theory to scattering of waves (electromagnetic, acoustic, etc.) in inhomogeneous
media. To be more precise, we consider in \S 4 the
operators 
\begin{equation} 
 H_0=M_0^{-1} (x) P(D)\quad
\mathrm{and} \quad H =M^{-1} (x) P(D).
\label{eq:perm}\end{equation} 
 Here $M_0 (x)$ and $M  (x)$ are
positively definite bounded matrix-valued functions and $P(D)$ is an
 elliptic differential operator. The operators $H_0$ and $H$ are self-adjoint in
Hilbert spaces with scalar products defined naturally in terms of $M_0$ and $M$, respectively.
Our main result is that   the wave operators for the pair $H_0$, $H$ exist, are isometric and
  complete  if
\begin{equation} 
|M(x)-M_0(x)|\leq C(1+|x|)^{-\rho}, \quad \rho >d, \quad x\in{\Bbb R}^d, 
\label{eq:VMMV}\end{equation}
($C$ is some constant). 
We emphasize that no special assumptions  on $M_0(x)$ are required,
that is the ``background" medium might be inhomogeneous. In contrast, the
smooth theory relies on a sufficiently explicit diagonalization of the operator $H_0$; for example, if
 $M_0$ does not depend on $x$, then $H_0$ can be diagonalized by the Fourier transform.
 On the other hand, in this
approach it suffices (see \cite{De}) to suppose that $\rho>1$ in the estimate (\ref{eq:VMMV}).

In \S 5 we study the scattering theory for the wave equation. This problem can be almost reduced
to that considered in \S 4. However the corresponding operator $P(D)$ is not a differential
operator and, if considered as a pseudo-differential operator, it has a non-smooth symbol. It
creates some new difficulties.

Both these problems were considered by  M. Sh. Birman in the papers \cite{B4a,B5}  
where  it was supposed that $M_0 (x)$ is a constant matrix (there was also a similar assumption
in the case of the wave equation). Essentially the same assumptions were made in the papers
\cite{SchWi, ReSi, Dei} (see also the book \cite{RS}).
Similarly to
\cite{B4a,B5}, we proceed from the conditions for the existence and completeness of the wave 
operators  established earlier in the paper
\cite{BB} by A. L. Belopolskii and M. Sh. Birman. However a verification of these conditions for
the pair (\ref{eq:perm}) in the case where $M_0$ is a function of $x$ is a relatively tricky
business. One of possible tricks is presented in this paper.

 Typically, one has to check that the operators
$\langle x \rangle^{-r}   (H_0-z)^{-n}$ and $\langle x \rangle^{-r}   (H -z)^{-n}$
belong  to the Hilbert-Schmidt class ${\goth S}_2$ if $r> d/2$ and $n$ is sufficiently large.
Since properties of the operators $H_0$ and $H$ are the same, we discuss only the
second of these operators. 
The problem here is that the inclusion
$\langle x \rangle^{-r}   (P(D)-z)^{-n}\in {\goth S}_p$ (${\goth S}_p$ is the Neumann-Schatten
class) implies that $\langle x \rangle^{-r}   (H-z)^{-n}\in {\goth S}_p$ for $n=1$ only.
If $n>1$, then a direct verification of this assertion requires restrictive assumptions
on derivatives of $M (x)$. Roughly speaking, we fix this difficulty in the following way.
We consider the whole scale of classes ${\goth S}_p$ and find such a number $p(r,n)$ that  $\langle x \rangle^{-r}  
(H-z)^{-n}\in {\goth S}_p$ for   $p>p(r,n)$. This is done successively for $n=1,2,\ldots$.
To make a passage from $n$ to $n+1$, we use that the operators
\[
(P(D)-z) \langle x \rangle^{-r}  (P(D)-z)^{-1} \langle x \rangle^r 
\] 
are bounded for all $r\geq 0$.  This allows us to deduce that
  $\langle x \rangle^{-r}  
(H-z)^{-n-1}\in {\goth S}_p$ for   $p>p(r,n)$ from the inclusions 
  $\langle x \rangle^{-r_0}  
(H-z)^{ -1}\in {\goth S}_p$ for   $p>p(r_0,1)$
and
  $\langle x \rangle^{-r_1}  
(H-z)^{-n }\in {\goth S}_p$ for   $p>p(r_1,n)$
with suitably chosen $r_0+r_1=r$.

 \section{Preliminaries}

{\bf 1}.
Let ${\cal H}_0$ and  ${\cal H}$ be two Hilbert spaces, and let ${\goth B}$ be the algebra of all bounded operators acting from
${\cal H}_0$ to  ${\cal H}$. The ideal of compact operators will be denoted by ${\goth S}_\infty$.
For any compact $A$ we denote by $s_n(A)$ the eigenvalues of the positive compact operator $(A^\ast A)^{1/2}$
listed with account of multiplicity in decreasing order.
Important
symmetrically normed ideals ${\goth S}_p$, $1\leq p<\infty$, of the algebra   ${\goth B}$ are
formed by operators
$A\in{\goth S}_\infty$ for which  
\[
 \sum_{n=1}^\infty s_n^p(A) <\infty.
\]
 In particular, ${\goth S}_1$ and ${\goth S}_2$ and called the trace and Hilbert-Schmidt  
classes, respectively.
  Clearly, ${\goth S}_{p_1}\subset {\goth S}_{p_2}$ for $p_1\leq p_2$.
 Moreover, we have 

\begin{proposition}\label{SS1}
 If  $A_j\in {\goth S}_{p_j}$, $j=1,2$, and $p^{-1}=p_1^{-1}+ p_2^{-1}\leq 1$, then $A=A_1A_2\in
{\goth S}_{p}$.
\end{proposition}

\medskip

{\bf 2.}
We need to consider integral operators of the form
\begin{equation}
(Tf)(x) =(2\pi)^{-d/2}\int_{{\Bbb R}^d} a(x)
  \exp( i\langle x,\xi\rangle) b(\xi) \hat{f}(\xi) d\xi 
\label{eq:IntOp}\end{equation}
 acting in the space $L_2 ({\Bbb R}^d ;{\Bbb C}^k)$.
Here 
\[ 
\hat{f}(\xi) =   (2\pi)^{-d/2}\int_{{\Bbb R}^d} \exp(-i\langle x,\xi\rangle)
f(x) dx
\]
is the Fourier transform   of $f\in L_2 ({\Bbb R}^d ;{\Bbb C}^k)$
and $ a(x)$, $b(\xi)$ are $k\times k$-matrix-functions which we always suppose to be bounded.
Then operator  (\ref{eq:IntOp}) is bounded. Below
we sometimes use the short-hand notation
$T=a(x)b(\xi)$ for operators of the form (\ref{eq:IntOp}). Let us also set
\[
\langle x \rangle =(1+ |x|^2)^{1/2}, \quad
\langle \xi\rangle =(1+ |\xi|^2)^{1/2}.
\]

The following assertion is well-known.

\begin{proposition}\label{4.1.5sim} 
The  operator $(\ref{eq:IntOp})$ is compact if the functions $a$ and $b$  tend to
zero at infinity. This operator
  belongs to the   class ${\goth S}_p(L_2({\Bbb R}^d ;{\Bbb C}^k))$, $p\geq 1$, if
\[
|a(x)|\leq C (1+|x|)^{-r},\quad 
|b(\xi)|\leq C (1+|\xi|)^{-r},\quad r > d/p.
\]
 \end{proposition}

 The proof of this result can be found, e.g., in \cite{RS}. Strictly speaking, the case
$p\in(1,2)$ was not considered in  \cite{RS}. However it can be directly deduced from
Proposition~\ref{4.1.5sim} for $p=1$ and $p=2$ by the complex interpolation.

\medskip

{\bf 3.}
 We need also conditions of boundedness   in the space $L_2({\Bbb R}^d ; {\Bbb C}^k )$ of products
of multiplication operators in $x$- and
$\xi$-representations.  

\begin{proposition}\label{1.9}
 Suppose that   a matrix-function $a(x)$ has
$n$ bounded derivatives and $0\leq l \leq n$. Then the product $\langle \xi\rangle^l a(x) \langle \xi\rangle^{-l}$
defined by its sesquilinear form on the Schwartz class
${\cal S}$ is a bounded operator. Its norm is estimated by
\[
\sup_{  |\sigma|\leq n}\sup_{x\in{\Bbb R}^d }  |(\partial^\sigma a) (x) |.
\]
\end{proposition}

\begin{pf} 
 Note first that, for all $j=1,\ldots, d$,
\begin{equation} 
 D_j^n a(x) \langle \xi\rangle^{-n} =\sum_{m=0}^n i^{-m} C_n^m \partial^m a (x)/\partial x_j^m
(\xi_j^{n-m} \langle\xi\rangle^{-n}),
\label{eq:1.30}\end{equation}
where $C_n^m$ are binomial coefficients. Since the functions 
$\partial^m a (x)/\partial x_j^m$  and $\xi_j^{n-m} \langle\xi\rangle^{-n}$
are bounded,  operator (\ref{eq:1.30}) is also bounded. This entails that
the operators  
$ |\xi_j|^n a(x) \langle \xi\rangle^{-n}$ and hence $\langle \xi\rangle^n a(x) \langle \xi\rangle^{-n}$
are bounded.

 To pass to an arbitrary $l$, we consider the
function
\[
 (a(x) \langle \xi\rangle^{-z} f, \langle \xi\rangle^z g), \quad f,g\in {\cal S},
\] 
analytic in $z$ and bounded in any strip $c_0\leq \Re z \leq c_1$. As we have
seen, this function is bounded by $C ||f||\: ||g||$ if $\Re z=n$ and of course if $\Re z=0$. In view of
the Hadamard three lines theorem (see, e.g., \cite{GK})  this implies that the same bound is true
for $z=l$.
\end{pf}

 Consider now a more general operator
\begin{equation} 
  \langle \xi\rangle^l a(x)\langle x\rangle^{-r} b(\xi)\langle \xi\rangle^{-l}
\langle x\rangle^r.
\label{eq:1.24}\end{equation}

\begin{proposition}\label{1.10}
 Suppose that matrix-functions  $a(x)$ and $b(\xi)$ have $n$ bounded
 derivatives and $0\leq  l\leq n$, $0\leq  r\leq n$.
 Then operator $(\ref{eq:1.24})$ 
defined by its sesquilinear form on the Schwartz class
${\cal S}$ is   bounded.
\end{proposition}

\begin{pf}
Set $\tilde{a}(x)= a(x) \langle x\rangle^{-r}$.
 Similarly to the
proof of Proposition~\ref{1.9}, we  consider the function
\[
 ( \tilde{a}(x) b(\xi)\langle \xi\rangle^{-z} \langle x\rangle^r f, \langle \xi\rangle^z g),
 \quad f,g\in {\cal S},
\]
 analytic in $z$ and bounded in any strip $c_0\leq \Re z \leq c_1$.    In view of
the Hadamard three lines theorem it suffices to verify that this function is
bounded by $C ||f||\: ||g||$ for $\Re z=n$ and for $\Re z=0$.
 According to (\ref{eq:1.30}) the operator
$   \xi_j^n \tilde{a}(x)  b(\xi)\langle \xi\rangle^{-n-i\alpha}\langle x\rangle^r$,
 $\alpha\in{\Bbb R}$,
is a sum of terms 
\begin{equation}
 (\partial^m \tilde{a} (x)/\partial x_j^m \langle x\rangle^r )\cdot (\langle x\rangle^{-r}    (
b(\xi)
\xi_j^{n-m}\langle
\xi\rangle^{-n-i\alpha})\langle x\rangle^r)
\label{eq:rescommm}\end{equation}
 where $m=0,1,\ldots, n$. The first factor here is a bounded function of $x$.
The function $b(\xi) \xi_j^{n-m}\langle \xi\rangle^{-n-i\alpha}$
is bounded, together with its $n$ derivatives, uniformly in $\alpha$. Therefore,  applying
Proposition~\ref{1.9} with the roles of the variables $x$ and $\xi$  interchanged
to the second factor in (\ref{eq:rescommm}),
 we see that  this operator is bounded uniformly in $\alpha$.
Similarly, the operators 
$   \langle x\rangle^{-r} b(\xi)\langle \xi\rangle^{-i\alpha}
\langle x\rangle^r$ are also bounded uniformly in $\alpha$.
\end{pf}

\medskip  

{\bf 4}.
let us consider self-adjoint operators $H_0$ and $H$ in Hilbert spaces
${\cal H}_0$ and ${\cal H} $, respectively.
Recall that the essential spectrum $\sigma^{ess} $ of $H$ is defined as its spectrum $\sigma$
without isolated eigenvalues of finite multiplicity.
The same objects for the operator $H_0$ will be always labelled by the index `$0$'. 
According to the Weyl theorem
$\sigma^{ess} =
\sigma^{ess}_0$ if ${\cal H}_0={\cal H} $ and the difference $H-H_0$ is a compact operator. We
note a simple generalization
 of this result. Below we use the notation $R_0(z)=(H_0-z)^{-1}$, $R (z)=(H -z)^{-1}$.

\begin{proposition}\label{Wey1}
Let an operator $J:{\cal H}_0  \rightarrow {\cal H}$ be bounded, and let the inverse operator $J^{-1}$
exist  and be also bounded.
 Suppose that   
\begin{equation}
R (z) J - JR_0 (z) \in {\goth S}_\infty   
\label{eq:rescomp1}\end{equation}
for some point $z\not\in \sigma_0 \cup\sigma  $. Then
$\sigma^{ess} =\sigma^{ess}_0$.
\end{proposition}

\begin{pf}
If $\lambda \in \sigma^{ess}_0$, then there exists a sequence
(Weyl sequence) $f_n$ such that $H_0 f_n-\lambda f_n \rightarrow 0$,
 $  f_n \rightarrow 0$  weakly as  $ n \rightarrow \infty$ and $||f_n || \geq c>0$.
 This ensures that
$R_0(z)f_n- (\lambda-z)^{-1} f_n \rightarrow 0$ and hence
$J R_0(z)f_n- (\lambda-z)^{-1} J f_n \rightarrow 0$.
Now it follows from condition  (\ref{eq:rescomp1}) that
\begin{equation}
R (z) Jf_n- (\lambda-z)^{-1} J f_n \rightarrow 0.
\label{eq:rescomp1xyz}\end{equation}
 Set   $g_n=R (z) J f_n$. Relation    (\ref{eq:rescomp1xyz}) means that
 $Hg_n-\lambda g_n \rightarrow 0$. Moreover, $g_n \rightarrow 0$ weakly as  $ n \rightarrow \infty$ and,
  by virtue of (\ref{eq:rescomp1xyz}), the relation $||g_n || \rightarrow 0$ would have implied that
$|| Jf_n || \rightarrow 0$ and therefore $||f_n || \rightarrow 0$. Thus, $g_n$ is a Weyl sequence
for the operator $H$ and the point
$\lambda$ so that $\sigma^{ess}_0 \subset \sigma^{ess }$.

To prove the opposite inclusion, we remark that
\[
J^{-1}R (z)  -  R_0 (z) J^{-1}\in {\goth S}_\infty 
\]
and use the result already obtained with the roles of $H_0$, $H$ interchanged and 
$J^{-1}$ in place of $J$.
\end{pf}

\section{Scattering with two Hilbert spaces}

{\bf 1.}
  Scattering theory requires classification of the spectrum in terms of the theory of
measure. Let $H$ be an arbitrary self-adjoint operator   in a Hilbert space $\Hilb$. We denote by
$E(\cdot)$ its spectral measure. Recall that 
there is a decomposition  $\Hilb =
\Hilb^{ac}\oplus\Hilb^{sc} \oplus\Hilb^{pp}$ into the orthogonal sum of invariant subspaces of the
operator $H$ such that the measures $(E(\cdot) f, f)$ are absolutely continuous, singular continuous
or pure point for all $f\in \Hilb^{ac}$,  $f\in \Hilb^{sc}$ or  $f\in \Hilb^{pp}$, respectively.
The operator $H$ restricted to $\Hilb^{ac}$,
$\Hilb^{sc}$ or $\Hilb^{pp}$ shall be denoted $H^{ac}$, $H^{sc}$ or $H^{pp}$, respectively. 
 The pure point   part corresponds to eigenvalues. The singular continuous part is typically
absent.  Scattering theory   studies the  absolutely continuous part $H^{ac}$ of $H$.
We denote $P $   the orthogonal projection onto the absolutely continuous subspace $\Hilb^{ac}$.

Let us consider the large time
behaviour of solutions
$$ 
u(t) = e^{-iH t} f.
$$
 of the time-dependent equation
$$
i \frac{ \partial u}{\partial t} = H u, \quad u(0) = f \in \Hilb.
$$
If $f$ is an eigenvector, $H f = \lambda f$, then   $u(t) = e^{-i \lambda t} f$, so the   time
behaviour is evident. By contrast, if $f \in \Hilb^{ac}$, one cannot, in general, calculate $u(t)$
explicitly, but scattering  theory allows us to find its asymptotics as $t\rightarrow\pm\infty$.
 In the perturbation  theory setting,
it is natural to understand the asymptotics of $u$ in terms of solutions of the unperturbed equation,
$i u_t = H_0 u$. To compare the operators $H_0$ and $H$, one has to introduce an `identification'
operator
$J:\Hilb_0 \rightarrow\Hilb$ which we suppose to be bounded. Suppose also that,
in some sense, $J$ is close to a unitary operator and the perturbation
 $HJ-J H_0$ is   `small'. Then it turns out that  for
all
$f\in\Hilb^{ac}$, there are
$f_0^\pm\in\Hilb_0^{ac}$ such that
\begin{equation}
\lim_{t \to \pm \infty}\big\| e^{-i H t} f - J e^{-i H_0 t} f_0^\pm \big\| = 0.
\label{WOP}\end{equation}
Hence   $f_0^\pm$ and $f$ are related by the equality   
$$
f = \lim_{t \to \pm \infty} e^{i H t} J e^{-i H_0 t} f_0^\pm,
$$
which justifies the following fundamental definition.   It goes back to C. M\o ller
for $\Hilb_0=\Hilb$ and $J=I$. In the general case it was formulated by
T. Kato \cite{Ka3}.

 \begin{definition} 
Let $J$ be a bounded operator. Then the   wave operator 
$W_\pm(H, H_0,\\ J)$ is defined by
\begin{equation}
 W_\pm(H, H_0 ; J) = \slim_{t \to \pm \infty} e^{i H t} J e^{-i H_0 t}P_0,
\label{WOJ}\end{equation}
 when this limit exists.
 \end{definition}
 
Clearly, relation (\ref{WOP}) holds for all $f$ from the range  $\Ran W_\pm$ of the wave operator
$W_\pm= W_\pm(H, H_0;J)$.  The wave operators enjoy the intertwining property  
$$ 
W_\pm(H, H_0; J) H_0 = H W_\pm(H, H_0; J).
$$ 
In our applications they are isometric  on $\Hilb_0^{ac}$ which is guaranteered by the condition
\[
  \slim_{t \to \pm \infty} (J^\ast J-I) e^{-i H_0 t}P_0=0.
\]
This implies that $H_0^{ac}$ is unitarily equivalent, via $W_\pm $, to the
restriction   of $H$ on the range $\Ran W_\pm$ of the wave operator $W_\pm$ and hence
$\Ran W_\pm\subset\Hilb^{ac}$. 

 \begin{definition}
Suppose that the wave operator $W_\pm(H, H_0; J)$ exists and is
isometric  on $\Hilb_0^{ac}$.
It is said to be complete if
\[
\Ran  W_\pm (H,H_0;J)={\cal H}^{ac}.
\]
 \end{definition}

 Thus, if $W_\pm(H, H_0, J)$ exists, is isometric and complete,
then $H_0^{ac}$ and $H^{ac}$ are unitarily equivalent.   

Suppose additionally that $J$ is boundedly invertible.
It is a simple result that $W_\pm(H, H_0, J)$ is complete if and only if the `inverse' wave
operator $W_\pm(H_0, H; J^{-1})$ exists.

\medskip

{\bf 2.}
Let us now discuss conditions for the existence of wave operators.
In the abstract framework they are given by the trace-class theory.
Its fundamental result  is the following theorem of Kato-Rosenblum-Pearson.

\begin{theorem}\label{Pe} 
Suppose that $H_0$ and $H$ are selfadjoint operators in spaces ${\cal
H}_0$ and ${\cal H}$, respectively,
$J:{\cal H}_0 \rightarrow{\cal H}$ is a bounded operator  and $V=HJ-JH_0\in{\goth S}_1$. Then the
WO $W_\pm (H,H_0;J)$ exist.
\end{theorem}

Actually, this result was established by T. Kato and M. Rosenblum  for the case
${\cal H}_0={\cal H}$, $J=I$ and then extended by D. Pearson
to the operators acting in different spaces.

Applications to differential operators require generalizations of this result. 
The following result of \cite{BB} gives
efficient conditions guaranteeing the existence of wave operators
and all their properties discussed above.
Its simplified proof relying on Theorem~\ref{Pe}  can be found in \cite{I}.

\begin{theorem}\label{6.6.8}
Suppose that the operator $J:{\cal H}_{0}\rightarrow {\cal H}$
has a bounded inverse and
$J{\cal D}(H_0)={\cal D}(H )$. 
Suppose that  
\begin{equation}
E(\Lambda)(HJ-JH_0) E_0(\Lambda) \in {\goth S}_1
\label{eq:JI1}\end{equation}
and
\begin{equation}
(J^\ast J -I) E_0(\Lambda) \in {\goth S}_\infty  
\label{eq:JI2}\end{equation}
for any bounded interval $\Lambda$. Then the WO
$W_\pm (H,H_0;J)$ exist, are isometric on ${\cal H}_{0}^{ac}$, and are complete.
Moreover, there exist the WO $W_\pm (H_0,H;J^\ast)$ and $W_\pm (H_0,H;J^{-1})$; these WO are equal
to one another  and to the operator $W_\pm^\ast (H,H_0;\\J)$; they are isometric on ${\cal H}
^{ac}$  and
  complete.
 \end{theorem}

 Thus, the absolutely continuous part of a self-adjoint
operator  is stable under fairly general perturbations. However assumptions on perturbations are
much more restrictive than those required for stability of the essential spectrum
(cf. Proposition~\ref{Wey1}).  

\section {Scattering problems for perturbations of a medium}

By definition, a matrix pseudodifferential operator  with constant coefficients acts in the momentum
representation as multiplication by some matrix-function. This function is called symbol  of such
pseudodifferential operator.

\medskip

  {\bf 1.}
Let $  L_2({\Bbb R}^d; {\Bbb C}^k)$ and let $P(D)=\Phi^\ast A\Phi$ where $A$ is multiplication
by a a symmetric $k\times k$- matrix-function
$A(\xi)$.
 The operator $A$ is
selfadjoint on domain ${\cal D}(A)$ which consists of functions
$\hat{f}\in L_2({\Bbb R}^d; {\Bbb C}^k)$ such that $A\hat{f}\in L_2({\Bbb R}^d; {\Bbb C}^k)$.  Hence
the operator $P(D)$ is selfadjoint on domain ${\cal D}(P(D))= \Phi^\ast{\cal D}(A)$. 
Below we need to restrict the
class of operators $P(D)$. Set
 \[
\nu(\xi)=\min_{|{\bold n}|=1}|A(\xi){\bold n}|,\quad {\bold n}\in {\Bbb C}^k. 
\]
Clearly, $\nu(\xi)$ is the smallest of the absolute values of eigenvalues of the matrix
$A(\xi)$. 
The operator $P(D)$ is called strongly Carleman if 
$\nu(\xi)\rightarrow\infty$ as
$|\xi|\rightarrow\infty$.
 We often need a stronger condition
\begin{equation}
\nu(\xi)\geq c |\xi|^\varkappa,\quad \varkappa>0,\quad c>0,
\label{eq:nu}\end{equation}
for $|\xi|$   sufficiently large.

 Clearly, $P(D)$ is a differential operator if entries of  the matrix $A(\xi)$
are polynomials of $\xi$. Let us denote by $A_0(\xi)$ the principal symbol of $A(\xi)$,
 that is $A_0(\xi)$ consists of terms of higher order which will be denoted by $\mathrm{ord}\,
P(D)$. If $\det A_0(\xi) \neq 0$ for $ \xi \neq 0$, then $P(D)$ is called
 elliptic   of order $\mathrm{ord}\, P(D)$.
For elliptic operators condition (\ref{eq:nu}) is satisfied with 
$\varkappa=\mathrm{ord}\, P(D)$.

Let $M_0(x)$ and $M(x)$ be    symmetric $k\times k$ - matrices satisfying the
condition
\begin{equation}
 0<c_0\leq M_0(x)\leq c_1<\infty, \quad 0<c_0\leq M (x)\leq c_1<\infty,
\label{eq:3.4.1}\end{equation}
and let $M_0$ and $M$ be the operators of multiplication by these matrices. 
We denote by  ${\cal
H} $ the Hilbert space with scalar product
\begin{equation}
 (f,g)_{\cal H} =\int_{{\Bbb R}^d}\langle M (x)f(x),g(x)\rangle_{{\Bbb C}^k} dx.
\label{eq:3.4.1M}\end{equation}
The space  ${\cal H}_0$ is defined quite similarly with $M (x)$ replaced
by $M_0(x)$.
 Of course the spaces ${\cal H}_{00}= L_2({\Bbb R}^d; {\Bbb C}^k)$, ${\cal H}_0$ and ${\cal H}$
consist of the same elements.
The operators $M_0$ and $M$ can   be considered in all these spaces.
The   operators $H_0$ and $H$ are defined  by the equalities  
(\ref{eq:perm})
on common domain  $ {\cal D}(H_0)={\cal D}(H )={\cal D} (P(D))$ in the spaces
${\cal H}_0$ and  ${\cal H}$, respectively. 
Their selfadjointness follows from selfadjointness of the operator $P(D)$ in the space
  $L_2({\Bbb R}^d; {\Bbb C}^k)$. 
 Let $I_0:{\cal H}_0\rightarrow{\cal H}$
and $I_1=I_0^{-1}:{\cal H} \rightarrow{\cal H}_0$ be the identical mappings.
They are
often omitted if this does not lead to any confusion.  Note however that
\begin{equation} 
I_0^\ast =M_0^{-1} M, \quad I_1^\ast =  M^{-1}M_0.
\label{eq:Iast}\end{equation}
Put $H_{00}=P(D)$, $R_{00}(z)=(H_{00}-z)^{-1}$.
Below we use the  resolvent identities
\begin{eqnarray}
 R (z) &= & R_{00}(z)
(M +z(M -I) R (z)), \quad z\not\in \sigma_{00}  \cup \sigma,
\label{eq:pm4}\\ 
 R (z) &= &   (I- z R(z)(M^{-1}-I)  ) R_{00}(z) M, \quad z\not\in \sigma_{00}  \cup \sigma,
\label{eq:pm4bi}\end{eqnarray}
 which can be verified by a direct calculation.
Of course   similar identities hold for $R_0(z)$. 

As far as the essential spectrum is concerned, we have the following standard assertion.

\begin{proposition}\label{PM1we}
Suppose that $H_{00}$ is strongly Carleman and 
\[ 
V(x):=M(x)-M_0(x)\rightarrow 0 
\]
as $| x | \rightarrow \infty$. 
  Then $\sigma^{ ess  } = \sigma^{ess}_0 $.
 \end{proposition}

\begin{pf}
Let us use the resolvent identity for the pair $H_0,H$:
\begin{equation}
 R (z) - R_0(z)
= R (z) M^{-1} V (I+ z R_0(z)), \quad z\not\in \sigma_0  \cup \sigma.
\label{eq:pm4bis}\end{equation}
 According to   identity (\ref{eq:pm4bi}) and Proposition~\ref{4.1.5sim},
the operators $R (z) M^{-1} V$ and hence   (\ref{eq:pm4bis}) are compact.
 Thus it remains to refer to
  Proposition~\ref{Wey1}. 
\end{pf}

\medskip

  {\bf 2.}
Let us pass to scattering theory.  
We proced from the following analytical result.

\begin{proposition}\label{PM3}
   Let $P(D)$ be an elliptic differential operator of order $\varkappa$ in the space ${\cal
H}=L_2({\Bbb R}^d; {\Bbb C}^k)$.   
Suppose that the function $M (x)$ obeys condition
$(\ref{eq:3.4.1})$. Set $H=M^{-1}P(D)$.
  Then  the operator $\langle
x\rangle^{-r}R^n (z)$,
$n=1,2,\ldots$, $z\not\in \sigma $, belongs to the class ${\goth S}_p$ provided   $p \geq 1$ and
\[
 p>d/\min\{r,\varkappa n\} =: p(r,n).
\]
 \end{proposition}

 \begin{pf}
The proof   proceeds by induction in $n$. If $n=1$, then we use the
equality
\begin{equation} 
\langle x\rangle^{-r}R   =(\langle x\rangle^{-r} \langle \xi\rangle^{- \varkappa})\cdot (\langle
\xi\rangle^{ \varkappa }R_{00} )\cdot((H_{00}-z) R ).
\label{eq:4xyz}\end{equation}
In the right-hand side the first factor belongs to the class ${\goth S}_p$ for 
 $ p > p (r,1)$   
 according to Proposition~\ref{4.1.5sim}. The second factor is a bounded operator
according to (\ref{eq:nu}), and the last factor is a bounded operator
according to (\ref{eq:pm4}).

To justify the passage from $n$ to $n+1$, we write the operator $\langle x\rangle^{-r}R^{n+1} $
as
\begin{eqnarray} 
 \langle x\rangle^{-r}R^{n+1} &=&(\langle x\rangle^{-r_0}R_{00})
\nonumber\\
&\times&
((H_{00}-z)\langle x\rangle^{-r_1}R_{00}\langle x\rangle^{r_1})\times 
(\langle x\rangle^{-r_1} (H_{00}-z)  R^{n+1} ),  
\label{eq:4.3.5}\end{eqnarray}
where $r_0 =r(n+1)^{-1}$, $r_1=n r_0$.  The first factor here 
belongs to  the class ${\goth S}_p$ for 
 $ p > p (r_0,1)$. The second factor 
 is bounded according to Proposition~\ref{1.10}.
It follows from   identity (\ref{eq:pm4}) that the last factor
\begin{equation}
\langle x\rangle^{-r_1}(H_{00}-z) R^{n+1} =\langle x\rangle^{-r_1} 
 M  R^n   + z \langle x\rangle^{-r_1} (M -I) R^{n+1}.
\label{eq:pm6}\end{equation}
This operator belongs to the class ${\goth S}_p$ where, by the inductive assumption, 
 $p> p(r_1, n)$.  Thus, by Proposition~\ref{SS1}, the product
(\ref{eq:4.3.5}) belongs to the class  ${\goth S}_p$ where 
\[
p^{-1}< p(r_0,1)^{-1} +p(r_1,n)^{-1}=(n+1) p(r_0,1)^{-1}=p(r, n+1)^{-1}
\]
 and of course $p\geq 1$.
\end{pf}

Now it is easy to prove

\begin{theorem}\label{PM4}
   Let $P(D)$ be an elliptic differential operator of order $\varkappa$ in the space ${\cal
H}=L_2({\Bbb R}^d; {\Bbb C}^k)$. 
 Assume that $M_0(x)$ and $M(x)$ satisfy conditions
$(\ref{eq:3.4.1})$ and $(\ref{eq:VMMV})$.    
  Then the  wave operators 
\[
W_\pm (H ,H_0; I_0), \quad W_\pm (H_0,H; I_0^\ast)
 \quad \mathrm{and} \quad W_\pm (H_0 ,H; I_1)
\]
 exist, are isometric
and are complete.
 \end{theorem}

\begin{pf}
By Theorem~\ref{6.6.8}, 
 we have only to check inclusions (\ref{eq:JI1}) and (\ref{eq:JI2}), that is  
\begin{equation}
E (\Lambda) (HI_0-I_0 H_0) E_0(\Lambda)\in{\goth S}_1  
\label{eq:pm3}\end{equation}
and
\begin{equation}
(I_0^\ast I_0 -I)E_0(\Lambda)\in{\goth S}_\infty 
\label{eq:pm31}\end{equation}
 for an arbitrary bounded interval $\Lambda$. It follows from (\ref{eq:perm}) and
(\ref{eq:Iast}) that    
\begin{equation}
H I_0-I_0 H_0=-M^{-1} V H_0,
\label{eq:perturb}\end{equation}
and
\begin{equation}
I_0^\ast I_0 -I= M_0^{-1} V.
\label{eq:isome}\end{equation}
By Proposition~\ref{PM3}, 
$V R_0^n \in{\goth S}_1$ if $n\varkappa > d$. Since  the operators 
$ H_0^n E_0(\Lambda)$ are bounded for all $n$, this implies both inclusions
(\ref{eq:pm3}) and (\ref{eq:pm31}) (actually, the operator in (\ref{eq:pm31}) also
 belongs to the trace class). 
\end{pf}

\begin{remark} 
Actually, we have verified that
$(HI_0-I_0 H_0) E_0(\Lambda)\in{\goth S}_1$
 which is stronger than (\ref{eq:pm3}). Inclusion (\ref{eq:pm3})
follows also from the inclusions
 $\langle x \rangle^{-r} E_0(\Lambda) \in{\goth S}_2$ and
 $\langle x \rangle^{-r} E (\Lambda) \in{\goth S}_2$ for $r> d/2$.
 \end{remark}

\section {Wave equation}

A propagation of sound waves  
   in inhomogeneous media is often described by the wave equation. Basically, the methods of the
previous section are applicable to this case. However, by a natural reduction of the wave equation
to the Schr\"odinger equation, the pseudodifferential operators with non-smooth symbols appear. This
requires a modification of Theorem~\ref{PM4}. Here we use the same notation as in the previous
section.

\medskip

  {\bf 1.}
Let us consider  the equation
\begin{equation} 
m(x)\frac{\partial^2 u(x,t)}{\partial t^2}=\Delta u(x,t),\quad x\in {\Bbb R}^d,
\label{eq:Weq}\end{equation}
where the function $m(x)$ satisfies condition
(\ref{eq:3.4.1}). Set
\begin{equation}  
 {\bold u}(x,t)=\left(\begin{array}{cc}   ((-\Delta)^{1/2}u)(x,t)  
 \\ \partial u(x,t)/ \partial t
\end{array}\right),\quad  M(x )=\left(\begin{array}{cc}  I & 0 
 \\ 0 & m(x)
\end{array}\right).
\label{eq:Weq1}\end{equation}
Then equation  (\ref{eq:Weq}) is equivalent to the   equation
\begin{equation} 
i M(x)\frac{\partial {\bold u}(x,t)}{\partial t }=
(-\Delta)^{1/2}\left(\begin{array}{cc}  0 & i  
 \\ -i  & 0
\end{array}\right)
{\bold u}(x,t).
\label{eq:Weq2}\end{equation}
According to (\ref{eq:Weq1}) initial data for equations
(\ref{eq:Weq}) and (\ref{eq:Weq2}) are connected by the relation
${\bold u}( 0)= (((-\Delta)^{1/2}u)(0), u_t(0))^t$
(the index $`t$' means `transposed').

 Set   
\begin{eqnarray*}  
 P(D)=(-\Delta)^{1/2}\left(\begin{array}{cc}  0 & i  
 \\ -i  & 0
\end{array}\right),
 \end{eqnarray*}
  and denote by  ${\cal
H} $ the Hilbert space with scalar product
(\ref{eq:3.4.1M}) where $k=2$. The operator $H =M^{-1}  P(D)$
is selfadjoint in the space  ${\cal H} $.
 Unitarity of the operator $\exp(-iHt)$
in this space   
is equivalent to the conservation of the energy
\[
||(-\Delta)^{1/2} u(t)ÊÊ||^2+ (m u_t (t), u_t (t)).
\]

Suppose now that another function   $m_0(x)$ also satisfying condition
(\ref{eq:3.4.1}) is given. All objects constructed by this function will be labelled by $`0$'.
Let $u_0(x,t)$ be a solution of equation (\ref{eq:Weq}) with $m (x)$ replaced by $m_0(x)$.
Our goal is to compare the asymptotics for large $t$ of solutions $u (x,t)$ and $u_0(x,t)$ in the
energy norm. This can be done in terms of the wave operators for the pair 
$H_0$, $H$. Indeed, we have the following obvious result.

\begin{proposition}\label{WEq}
Let $f = (((-\Delta)^{1/2}u)(0), u_t(0))^t$,
 $f_0 = (((-\Delta)^{1/2}u_0)(0), u_{0,t}(0))^t$
and let $t\rightarrow \infty$ $($or $t\rightarrow -\infty)$. Then relations
\[ 
||\exp(-iH t) f-I_0 \exp(-iH_0 t)f_0ÊÊ||_{\cal H}
\rightarrow 0 
\]
and
\[ 
||(-\Delta)^{1/2} ( u(t)-u_0(t)Ê)Ê||\rightarrow 0,
\quad
||    u_t(t)-u_{0,t}(t)Ê)Ê||\rightarrow 0
 \]
are equivalent to each other.
 \end{proposition} 
 
\medskip

  {\bf 2.}
According to Proposition~\ref{WEq} scattering theory for the wave equation reduces to a proof of
the existence and completeness of the wave operators
$W_\pm (H ,H_0; I_0)$. Now the symbol of the operator $P(D)$ equals
\[
A(\xi)=|\xi| \left(\begin{array}{cc}  0 & i  
 \\ -i  & 0
\end{array}\right).
\]
This function has a singularity at $\xi=0$ so that
 Theorem~\ref{PM4} cannot be directly applied. 
Since however this singularity is not too strong,  we have the following result.

\begin{theorem}\label{WEq1}
Let $d\leq 3$. Assume that   functions  $m_0(x)$ and $m(x)$ satisfy the condition
\[
 0<c_0\leq m_0(x)\leq c_1<\infty, \quad 0<c_0\leq m (x)\leq c_1<\infty 
\]
 and that  
\[
|m(x)-m_0(x)|\leq C (1+|x|)^{-\rho}, \quad \rho > d.
\]   
  Then all conclusions of  Theorem~$\ref{PM4}$ hold.
 \end{theorem}

 \begin{pf}
Again by Theorem~\ref{6.6.8}, 
 we have only to verify the inclusions   
(\ref{eq:pm3}) and (\ref{eq:pm31}). It follows from (\ref{eq:perturb}) and (\ref{eq:isome})
 that   it suffices to verify the inclusion
\begin{equation}
\langle x \rangle^{-r} E (\Lambda) \in {\goth S}_2,\quad r=\rho/2>d/2,  
\label{eq:Weq6}\end{equation}
and the same inclusion for the operator $H_0$. The operators
$H_0$ and $H$ are quite symmetric, and hence
we have to check (\ref{eq:Weq6}) only.

Let first $d=1$. The operator
$\langle x \rangle^{-r}\langle \xi \rangle^{-1}  \in {\goth S}_2$
so that, by virtue of (\ref{eq:4xyz}) where $\varkappa=1$, the operator
$\langle x \rangle^{-r} R   $  is
also Hilbert-Schmidt.

In the cases $d=2$ and $d=3$  we check that
\begin{equation}
\langle x \rangle^{-r} R ^2  \in {\goth S}_2. 
\label{eq:Weq7}\end{equation}
 Using again (\ref{eq:4.3.5}), (\ref{eq:pm6})  for $n=1$ and $r_0=r_1=r/2$, we see that
\begin{eqnarray*}
\langle x \rangle^{-r} R ^2&=&(\langle x \rangle^{-r/2} R_{00})\,
(( H_{00}-z)\langle x \rangle^{-r/2} R_{00} \langle x \rangle^{ r/2})
 \\
&\times&
(\langle x \rangle^{-r/2} (M + z (M  -I) R ) R ). 
\label{eq:Weq8}\end{eqnarray*}
By Proposition~\ref{4.1.5sim}, the operators $\langle x \rangle^{-r/2} R_{00}$
and hence $\langle x \rangle^{-r/2}R $ belong to the class  ${\goth S}_4$.
Therefore it remains to notice that the function
$(A(\xi)-z)^{-1} \langle \xi \rangle$
is bounded together with its first derivatives so that, by Proposition~\ref{1.10},
\[
( H_{00}-z)\langle x \rangle^{-r/2} R_{00} \langle x \rangle^{ r/2} 
 \in {\goth B}, \quad r\leq 2. 
\]
This yields (\ref{eq:Weq7}).
\end{pf}


\begin{thebibliography}{hhhh}



\bibitem {BB}A. L. Belopolskii and M. Sh. Birman,   The existence of   wave
operators in scattering theory in a couple of spaces,  Math. USSR Izv. {\bf 2}  (1968),
1117-1130.
 
 

\bibitem {B4a}M. Sh. Birman,  Some applications of a local condition 
 for the existence of   wave operators, Soviet Math. Dokl. {\bf 10} (1969), 393-397.

 \bibitem {B5}M. Sh. Birman,   Scattering problems for differential operators with perturbation
of the space, Math. USSR Izv. {\bf 5}  (1971), 459-474.

 

\bibitem {BY}M. Sh. Birman and D. R. Yafaev,  On the trace-class method in potential scattering
theory, J. Soviet Math. {\bf 56} no. 2 (1991), 2285-2299.


\bibitem {GK} I. C. Gokhberg and M. G. Kre\u{\i}n, {\em Introduction to the theory of linear
nonselfadjoint operators in Hilbert space},  Amer. Math. Soc., Providence, R. I., 1970. 


\bibitem {De}V. G. Deich,   The completeness of wave operators for systems with uniform propagation,
Zap. nauchn. sem. LOMI   {\bf 22} (1971), 36-46  (Russian). 


\bibitem {Dei}P. Deift,    Applications of a commutation formula, Duke Math. J. {\bf
45} (1978), 267-310.

\bibitem {Ka3}T. Kato,  Scattering theory with two Hilbert spaces, J. Funct. Anal. {\bf 1 }
(1967), 342-369. 
 
 

\bibitem {Pe}D. B. Pearson,  A generalization of the Birman trace theorem, J. Funct. Anal. {\bf
28} (1978), 182-186. 
 

\bibitem {ReSi}  M. Reed and B. Simon,  The scattering of classical waves from inhomogeneous
media,   Math. Z. {\bf 155} (1977), 163-180.

\bibitem {RS} M. Reed and B. Simon, {\em Methods of modern mathematical physics}, Vol   3, 
Academic Press, San Diego, CA,   1979.


\bibitem {SchWi} J. R. Schulenberger and C. H. Wilcox,   Completeness of the wave operators for
perturbations of uniformly propagative systems, J. Funct. Anal. {\bf 7} (1971), 447-474. 

\bibitem{I} D. R. Yafaev, {\em Mathematical scattering theory}, Amer. Math. Soc., Providence,
Rhode Island, 1992.
 

\bibitem{Y2} D. Yafaev, {\em Scattering theory: some old and new problems},
 Lecture Notes in Mathematics 1735, Springer, 2000.

 





 

\end{thebibliography}
\end{document}